\documentclass[reqno,12pt]{amsart}
\usepackage{amsmath,amsthm,amssymb,amsfonts,amscd}
\usepackage{xypic}
\theoremstyle{plain} 
	\newtheorem{thm}{Theorem}[section]
	\newtheorem*{thm*}{Theorem}
	
	\newtheorem{lem}[thm]{Lemma}
	
	\newtheorem{prop}[thm]{Proposition}
	\newtheorem{conj}[thm]{Conjecture}
	\newtheorem*{conj*}{Conjecture}
	
\theoremstyle{definition}
	\newtheorem{defn}[thm]{Definition}

\theoremstyle{remark}
	\newtheorem{rem}[thm]{Remark}
	
	\newtheorem*{pf}{Proof}
\numberwithin{equation}{section}
\def\CC{{\mathbb C}}

\def\LL{{\mathbb L}}

\def\QQ{{\mathbb Q}}
\def\RR{{\mathbb R}}
\def\ZZ{{\mathbb Z}}

\def\D{{\mathcal D}}

\def\K{{\mathcal K}}
\def\L{{\mathcal L}}

\def\N{{\mathcal N}}

\def\P{{\mathcal P}}

\def \mf#1#2#3#4{
\xymatrix{
{#1}\  \ar@<0.4ex>[r]^{{#2}} & \ {#4}
\ar@<0.4ex>[l]^{{#3}}
}
}

\def \mfs#1#2#3#4{\!
\xymatrix@C=1,5em{{#1} \! \ar@<0.2ex>[r]^{{#2}} & \! {#4}
\ar@<0.2ex>[l]^{{#3}}
}
\!}

\def \mfl#1#2#3#4{
\xymatrix@C=2.6em{{#1}\  \ar@<0.4ex>[r]^{{#2}} &\  {#4}
\ar@<0.2ex>[l]^{{#3}}
}
}

\def \mfss#1#2#3#4{\!
\xymatrix@C=1.5em{{#1} \ar@<0.3ex>[r]^{{#2}} & {#4}
\ar@<0.3ex>[l]^{{#3}}
}
\!}
\begin{document}
\title{On the categorical entropy and the topological entropy}
\date{\today}
\author{Kohei Kikuta}
\address{Department of Mathematics, Graduate School of Science, Osaka University, 
Toyonaka Osaka, 560-0043, Japan}
\email{k-kikuta@cr.math.sci.osaka-u.ac.jp}
\author{Atsushi Takahashi}
\address{Department of Mathematics, Graduate School of Science, Osaka University, 
Toyonaka Osaka, 560-0043, Japan}
\email{takahashi@math.sci.osaka-u.ac.jp}
\begin{abstract}
To an exact endofunctor of a triangulated category with a split-generator,
the notion of entropy is given by Dimitrov--Haiden--Katzarkov--Kontsevich, 
which is a (possibly negative infinite) real-valued function of a real variable. 
In this paper, we propose a conjecture which naturally generalizes the theorem by Gromov--Yomdin, and show that the categorical entropy of a surjective endomorphism of a smooth projective variety is equal to its topological entropy. 
Moreover, we compute the entropy of autoequivalences of the derived category in the case of the ample canonical or anti-canonical sheaf. 
\end{abstract}
\maketitle
\section{Introduction}
The topological entropy of a self-map of a manifold is the fundamental numerical invariant encoding
``complexity" of the map, which has been extensively studied by several people.
Recently, to an exact endofunctor of a triangulated category with a split-generator, 
the notion of entropy is given by Dimitrov--Haiden--Katzarkov--Kontsevich in \cite{DHKK}. 
It is a (possibly negative infinite) real-valued function of a real variable motivated by an analogy with the topological entropy. 

It is known that the topological entropy of a surjective holomorphic endomorphism of a compact K\"{a}hler manifold coincides with the natural logarithm of the spectral radius of the induced action on the cohomology, 
which is the fundamental theorem of Gromov-Yomdin \cite{Gro1,Gro2,Yom} (see also Theorem 3.6 in \cite{Ogu}). 
In particular, this enables us to calculate the topological entropy of an automorphism of a smooth projective variety over $\CC$ 
inside of algebraic geometry. 

The bounded derived category of coherent sheaves is a ``primitive" invariant of a smooth proper variety over $\CC$; 
there are many interesting examples of non-isomorphic varieties with an equivalent derived category.
Therefore, it is important whether the topological entropy is determined by the derived category or not. 
Concerning the entropy at zero of an exact endofunctor of the derived category, 
there is a lower bound, under a certain technical assumption, by the natural logarithm of the spectral radius of the induced action on the Hochschild homology (Theorem~2.9 in \cite{DHKK}).
Moreover, if the variety is projective then the categorical entropy and the topological entropy of a surjective endomorphism satisfying the assumption coincides (Theorem 2.12 in \cite{DHKK}). 

In this paper, we propose a conjecture which naturally generalizes the theorem by Gromov--Yomdin 
and we show the following two theorems:
\begin{thm*}[Theorem~\ref{main1}]
Let $X$ be a smooth projective variety over $\CC$ and $\D^b(X)$ the bounded derived category of coherent sheaves on $X$. 
Let $f:X\longrightarrow X$ be a surjective endomorphism 
and $\LL f^*:\D^b(X)\longrightarrow \D^b(X)$ the derived functor of $f$. 
Then, the categorical entropy $h_{cat}(f)$ coincides with the natural logarithm of the spectral radius $\rho([\LL f^*])$ of the induced automorphism $[\LL f^*]$ on the numerical Grothendieck group $\N(X)$ of $\D^b(X)$. Moreover, we have
\[
h_{cat}(f)=h_{top}(f),
\]
where $h_{top}(f)$ is the topological entropy of $f$.
\end{thm*}
\begin{thm*}[Theorem~\ref{main2}]
Let $X$ be a smooth projective variety over $\CC$ with the ample canonical or anti-canonical sheaf and 
$F$ an autoequivalence of the bounded derived category $\D^b(X)$ of coherent sheaves on $X$. We have
\[
h(F)=\log\rho([F])=0.
\]
\end{thm*}

The contents of this paper is as follows. 
In Section~2 and~3, we prepare notations and recall the definition and some properties of the entropy of exact endofunctors by \cite{DHKK}. 
After recalling the notion of dynamical degrees and the fundamental theorem due to Gromov--Yomdin in Section~4, 
we propose a conjecture (Conjecture~\ref{conj}) which naturally generalizes it and state our results in Section 5, 
whose proofs are given in Section~6. 

\bigskip
\noindent
{\it Acknowledgment}\\
\indent
The first and second named authors would like to thank Keiji Oguiso and Yuuki Shiraishi for valuable discussions.
The first named author is supported by JSPS KAKENHI Grant Number JP17J00227.
The second named author is supported by JSPS KAKENHI Grant Numbers JP24684005, JP26610008 and JP16H06337.
\section{Notations}
Throughout this paper,  we denote by $X$ a smooth projective variety over $\CC$ of dimension $d$ and 
by $\D^b(X)$ the bounded derived category of coherent sheaves on $X$ whose translation functor is given by $[1]$.
An endomorphism of $X$ is a morphism $f:X\longrightarrow X$ and an automorphism of $X$ is an endomorphism of $X$ 
which is invertible.
We denote the semi-group~(resp.~group) of endomorphisms~(resp.~automorphisms) of $X$ 
by ${\rm End}(X)$~(resp.~${\rm Aut}(X)$). 
An endofunctor of $\D^b(X)$ is a $\CC$-linear functor $F:\D^b(X)\longrightarrow \D^b(X)$ which is exact (or triangulated), 
namely, is a functor which commutes with the translation functor and sends exact triangles to exact triangles.
An endofunctor $F$ of $\D^b(X)$ is of {\it Fourier--Mukai type} 
if there exists an object $\P\in\D^b(X\times X)$ such that
$F(-)\cong \RR\pi_{2*}(\P\otimes_{X\times X}^\LL \LL\pi_1^*(-))$ 
where $\pi_1,\pi_2$ are the projections to the first and second components.

\section{Entropy of endofunctors}
It is well-known that $\D^b(X)$ is a split-closed triangulated category of finite type and it has a split-generator. 
Namely, it contains all direct summands of its objects, 
$\sum_{n\in\ZZ}\dim_\CC{\rm Hom}_{\D^b(X)}(M,N[n])<\infty$ for all $M,N\in\D^b(X)$ 
and there exists an object $G$ such that $\D^b(X)=\langle G\rangle$ where $\langle G\rangle$ is 
the smallest triangulated subcategory containing $G$ closed under isomorphisms and taking direct summands.
In particular, for a very ample invertible sheaf $\L$ on $X$, the objects $G:=\oplus_{i=1}^{d+1}\L^i$ and its dual $G^*:=\oplus_{i=1}^{d+1}\L^{-i}$ are split-generators (cf. Theorem 4 in \cite{Orl}). 
\begin{defn}[Definition~2.1 in \cite{DHKK}]
For each $M,N\in\D^b(X)$, define the function $\delta_{t}(M,N):\RR\longrightarrow\RR_{\geq0}\cup\{ \infty \}$ in $t$ by
{\small 
\begin{equation*}
\delta_{t}(M,N):= 
\begin{cases}
0
 & \text{ if }N\cong0\\
\inf\left\{
\displaystyle\sum_{i=1}^p {\rm exp}(n_i t)~
\middle
|~
\begin{xy}
(0,5) *{0}="0", (20,5)*{A_{1}}="1", (30,5)*{\dots}, (40,5)*{A_{p-1}}="k-1", (60,5)*{N\oplus N'}="k",
(10,-5)*{M[n_{1}]}="n1", (30,-5)*{\dots}, (50,-5)*{M[n_{p}]}="nk",
\ar "0"; "1"
\ar "1"; "n1"
\ar@{.>} "n1";"0"
\ar "k-1"; "k" 
\ar "k"; "nk"
\ar@{.>} "nk";"k-1"
\end{xy}
\right\}
 & \text{ if }N\in\langle M\rangle \\
\infty
 & \text{ if }N\not\in\langle M\rangle.
\end{cases}
\end{equation*}
}
The function $\delta_{t}(M,N)$ is called the {\it complexity} of $N$ with respect to $M$.
\end{defn}
\begin{rem}
If $M\in\D^b(X)$ is not isomorphic to a zero object, 
then an inequality $1\leq\delta_{0}(G,M)<\infty$ holds for any split-generator $G\in \D^b(X)$.
\end{rem} 
We recall some basic properties of the complexity.  
\begin{lem}[Proposition~2.3 in \cite{DHKK}]\label{complexity}
Let $M_1,M_2,M_3\in \D^b(X)$.
\begin{enumerate}
\item
If $M_1\cong M_3$, then $\delta_{t}(M_1,M_2)=\delta_{t}(M_3,M_2)$. 
\item
If $M_2\cong M_3$, then $\delta_{t}(M_1,M_2)=\delta_{t}(M_1,M_3)$. 
\item
If $M_2\not\cong 0$, then $\delta_{t}(M_1,M_3)\leq\delta_{t}(M_1,M_2)\delta_{t}(M_2,M_3)$.
\item
For any endofunctor $F$ of $\D^b(X)$, we have $\delta_{t}(F(M_1),F(M_2))\leq\delta_{t}(M_1,M_2)$.
\end{enumerate}
\end{lem}
\begin{defn}[Definition~2.5 in \cite{DHKK}]
Let $G$ be a split-generator of $\D^b(X)$ and $F$ an endofunctor of $\D^b(X)$ such that $F^n G\not\cong 0$ for all $n> 0$. 
The {\it entropy} of $F$ is the function $h_t(F):\RR\longrightarrow\{ -\infty \}\cup\RR$ given by
\begin{equation}
h_{t}(F):=\displaystyle\lim_{n\rightarrow\infty}\frac{1}{n}\log \delta_{t}(G,F^{n}G).
\end{equation}
\end{defn}
It follows from Lemma~2.6 in \cite{DHKK} that the entropy is well-defined.
\begin{rem}
If there exists a positive integer $n$ such that $F^n G\cong 0$, then $F^n M\cong 0$ for all $M\in\D^b(X)$.
\end{rem}
\begin{lem}\label{two-split-gen}
Let $G,G'$ be split-generators of $\D^b(X)$ and $F$ an endofunctor of $\D^b(X)$ such that $F^n G,F^nG'\not\cong 0$ for all $n> 0$.
Then
\begin{equation}
h_t(F)=\displaystyle\lim_{n\rightarrow\infty}\frac{1}{n}\log \delta_{t}(G,F^{n}G').
\end{equation}
\end{lem}
\begin{pf}
It follows from Lemma~\ref{complexity} (iii) and (iv) that
\[
\delta_{t}(G,F^nG')\leq \delta_{t}(G,F^nG)\delta_{t}(F^nG,F^nG')\leq \delta_{t}(G,G')\delta_{t}(G,F^nG).
\]
Similarly, we have $\delta_{t}(G,F^nG)\leq\delta_{t}(G',G)\delta_{t}(G,F^nG')$, which yields the statement.
\qed
\end{pf}
\begin{lem}[See Section~2 in \cite{DHKK}]\label{entropy}
Let $G$ be a split-generator of $\D^b(X)$ and $F_1,F_2$ endofunctors of $\D^b(X)$ 
such that $F_i^n G\not\cong 0$ for $i=1,2$ and for all $n> 0$. 
\begin{enumerate}
\item
If $F_{1} \cong F_{2}$, then $h_t(F_{1})=h_t(F_{2})$. 
\item
We have $h_t(F_{1}^m)= mh_t(F_{1})$ for $m\geq1$.
\item
We have $h_t([m])=mt$ for $m\in\ZZ$. 
\item
If $F_1F_2\cong F_2F_1$, then $h_t(F_1F_2)\leq h_t(F_1)+h_t(F_2)$.
\item
If $F_1=F_2[m]$ for $m\in\ZZ$, then $h_t(F_1)=h_t(F_2)+mt$. 
\end{enumerate}
\end{lem}
The following proposition enables us to compute entropies. 
\begin{prop}[cf. Theorem~2.7 in \cite{DHKK}]\label{smooth-proper}
Let $G,G'$ be split-generators of $\D^b(X)$ and $F$ an endofunctor of $\D^b(X)$ of Fourier--Mukari type 
such that $F^n G,F^nG'\not\cong 0$ for all $n> 0$. 
Then the entropy $h_t(F)$ is given by
\begin{equation}
h_t(F)=\lim_{n\rightarrow\infty}\frac{1}{n}\log\delta'_{t}(G,F^nG'),
\end{equation}
where 
\begin{equation}
\delta'_{t}(M,N):=\sum_{m\in\ZZ}\left(\dim_{\CC} {\rm Hom}_{\D^b(X)}(M,N[m])\right)\cdot e^{-mt},\quad M,N\in\D^b(X).
\end{equation}
\end{prop}
\begin{pf}
The following is proven in the proof of Theorem~2.7 in \cite{DHKK}.
\begin{lem}\label{delta-delta'}
There exist $C_1(t),C_2(t)$ for $t\in\RR$ such that
\[
C_1(t)\delta_{t}(G,M)
\leq\delta'_{t}(G,M)
\leq C_2(t)\delta_{t}(G,M),\quad M\in\D^b(X).
\]
In particular, for each $M\in \D^b(X)$ we have
\begin{equation}
\lim_{n\rightarrow\infty}\frac{1}{n}\log\delta_{t}(G,M)
=\lim_{n\rightarrow\infty}\frac{1}{n}\log\delta'_{t}(G,M).
\end{equation}
\end{lem}
Together with Lemma~\ref{two-split-gen}, we have
\[
h_t(F)=\lim_{n\rightarrow\infty}\frac{1}{n}\log\delta_{t}(G,F^nG')
=\lim_{n\rightarrow\infty}\frac{1}{n}\log\delta'_{t}(G,F^nG').
\]
We finished the proof of the proposition.
\qed
\end{pf}
Proposition~\ref{smooth-proper} yields the following
\begin{lem}[Lemma~2.14 in \cite{DHKK}]\label{pic}
For each invertible sheaf $\L'\in {\rm Pic}(X)$, the entropy $h(-\otimes \L')$ of the autoequivalence functor $-\otimes \L'$ is zero. 
\end{lem}
\section{Topological entropy}
Let $X$ be a $d$-dimensional complex smooth projective variety, as in Section 2. 
Let $f$ be a surjective endomorphism of $X$. 
For an integer $q~(0\le q \le d)$, define $H^{q,q}(X;\RR):=H^{q,q}(X)\cap H^{2q}(X;\RR)$. 
Note that $f$ induces an automorphism $f^*$ of $H^{q,q}(X;\RR)$. 
\begin{defn}
For an integer $q~(0\le q \le d)$, define the {\it$q$-th dynamical degree} $d_q(f)$ of $f$ by
\begin{equation}
d_q(f):=\limsup_{n\to\infty}\left(\int_X (f^n)^* c_1(\L)^q\cup c_1(\L)^{d-q}\right)^{\frac{1}{n}}, 
\end{equation}
where $\L$ is a very ample invertible sheaf on $X$.   
\end{defn}
\begin{rem}
The $q$-th dynamical degree does not depend on the choice of $\L$.
\end{rem}
Decompose $H^{q,q}(X)=H^{q,q}(X;\RR)\otimes \CC$ into the direct sum of the $f^*$-invariant subspaces 
\begin{equation}\label{growth1}
H^{q,q}(X)=H^{q,q}(X)_{\lambda_1,m_1}\oplus\dots\oplus H^{q,q}(X)_{\lambda_s,m_s},\quad (|\lambda_1|,m_1)\ge \dots \ge (|\lambda_s|,m_s),
\end{equation}
where $f^*|_{H^{q,q}(X)_{\lambda_i,m_i}}$ is given by the Jordan block $J_{\lambda_i,m_i}$ and we write
$(|\lambda_i|,m_i)\ge (|\lambda_j|,m_j)$ if either $|\lambda_i|>|\lambda_j|$ or $|\lambda_i|=|\lambda_j|$ and $m_i\ge m_j$.
Note that for any matrix norm $\|\cdot\|$, we have
\begin{equation}\label{growth2}
\limsup_{n\to\infty}\frac{\|J_{\lambda_i,m_i}^n\|}{n^{m_i-1}\lambda_i^n}<\infty,\quad 
\limsup_{n\to\infty}\frac{n^{m_i-1}\lambda_i^n}{\|J_{\lambda_i,m_i}^n\|}<\infty.
\end{equation}
\begin{defn}[cf. Section~2.2 in \cite{DS}]
Let the notations as above.
The {\it spectral radius} $r_q(f)$ of the automorphism $f^*$ of $H^{q,q}(X;\RR)$ is the positive number $|\lambda_1|$, 
the maximum of absolute values of eigenvalues of $f^*$.
The {\it multiplicity} $l_q(f)$ of the spectral radius $r_q(f)$ is the integer $m_1$.
\end{defn}
\begin{prop}[cf. Proposition~2.5 in \cite{DS}]\label{growth-p}
Let $\L$ be a very ample invertible sheaf on $X$.   
For  an integer $q~(0\le q \le d)$, 
\begin{equation}
\limsup_{n\to\infty}\frac{\int_X (f^*)^n c_1(\L)^q\cup c_1(\L)^{d-q}}{n^{l_q(f)-1}r_q(f)^n}<\infty,\quad 
\limsup_{n\to\infty}\frac{n^{l_q(f)-1}r_q(f)^n}{\int_X (f^*)^n c_1(\L)^q\cup c_1(\L)^{d-q}}<\infty.
\end{equation}
In particular, $(\int_X (f^*)^n c_1(\L)^q\cup c_1(\L)^{d-q})^{1/n}$ converges to the $q$-th dynamical degree $d_q(f)$, 
which coincides with the spectral radius $r_q(f)$.
\end{prop}
\begin{rem}
By the Perron--Frobenius theorem, we can show the existence of $i~(1\le i\le s)$ such that $(\lambda_i,m_i)=(r_q(f),l_q(f))$. Moreover, since $f^*$ respects the subspace $H^{q,q}(X;\ZZ):=H^{q,q}(X)\cap H^{2q}(X;\ZZ)$ (which is non-empty since $X$ is projective), 
it turns out that $r_q(f)$ is an algebraic integer. 
\end{rem}
Let $r(f)$ be the {\it spectral radius} of the automorphism $f^*$ of the (total) cohomology group $H^\bullet(X;\CC)$, 
the maximum of absolute values of eigenvalues of $f^*$.
The following is the fundamental theorem on topological entropy due to Gromov-Yomdin \cite{Gro1,Gro2,Yom} (see also Theorem 3.6 in \cite{Ogu}). 
\begin{prop}\label{gy}
For each surjective endomorphism $f\in {\rm End}(X)$,
\begin{equation}
h_{top}(f)=\log r(f)=\log\max_q r_q(f)=\log\max_q d_q(f).
\end{equation}
\end{prop}
\begin{rem}
Since $f^*$ is an automorphism of $H^{q,q}(X;\RR)$ and $f^*(H^{q,q}(X;\ZZ))\subset H^{q,q}(X;\ZZ)$, 
we have $r_q(f)\ge 1$ for all $q~(0\le q \le d)$ and hence $h_{top}(f)\ge 0$. 
\end{rem}
\section{Results}
\subsection{Spectral radii}
In order to study the structure of the entropy of endofunctors on $\D^b(X)$ we prepare some terminologies.
For $M,N\in\D^b(X)$, set 
\begin{equation}
\chi(M,N):=\sum_{n\in\ZZ}(-1)^n\dim_\CC{\rm Hom}_{\D^b(X)}(M,N[n]).
\end{equation}
It naturally induces a bilinear form on the Grothendieck group $K_0(X)$ of $\D^b(X)$, called the {\it Euler form}, 
which is denoted by the same letter $\chi$. 
Then {\it numerical Grothendieck group} $\N(X)$ is defined as the quotient of $K_0(X)$ by the radical of $\chi$ 
(which is well-defined by the Serre duality). 
It is important to note that $\N(X)$ is a free abelian group of finite rank by Hirzebruch-Riemann-Roch theorem 
and that $\N(X)_\QQ:=\N(X)\otimes_\ZZ\QQ$ can be identified with 
$A_{num}^\bullet(X)_\QQ:=A_{num}^\bullet(X) \otimes_\ZZ\QQ$ 
where $A_{num}^q(X)$ is the Chow group of codimension $q$ cycles on $X$ modulo numerical equivalence. 

Let ${\rm End}^{FM}(\D^b(X))$ be the semi-group of isomorphism classes of endofunctors of $\D^b(X)$ of Fourier--Mukai type 
and ${\rm End}(\N(X))$ the semi-group of endomorphisms of $\N(X)$. 
Since an endofunctor of Fourier--Mukai type has left and right adjoint functors, it respects the radical of $\chi$.
Therefore, we have the semi-group homomorphism
\begin{equation}
{\rm End}^{FM}(\D^b(X))\rightarrow{\rm End}(\N(X)),\quad F\mapsto [F].
\end{equation}
Define $\N(X)_\CC:=\N(X)\otimes_\ZZ\CC$ and denote by ${\rm End}_\CC(\N(X)_\CC)$ 
the semi-group of $\CC$-linear endomorphisms of $\N(X)_\CC$. 
For each $F\in {\rm End}^{FM}(\D^b(X))$, define $\rho([F])$ to be the {\it spectral radius} of $[F]$ 
considered as an element in ${\rm End}_\CC(\N(X)_\CC)$.

The {\it Hochschild homology group} of $X$ is defined as $H\!H_\bullet(X):=H\!H_\bullet(\D^b_{dg}(X))$ 
where $\D^b_{dg}(X)$ is a dg enhancement of $\D^b(X)$.  
We also have the semi-group homomorphism
\begin{equation}
{\rm End}^{FM}(\D^b(X))\rightarrow{\rm End}_\CC(H\!H_\bullet(X)),\quad F\mapsto H\!H_\bullet(F),
\end{equation}
where $H\!H_\bullet(F)$ is the induced $\CC$-linear endomorphism $H\!H_\bullet(F)$ of $H\!H_\bullet(X)$.
Denote by $\rho(H\!H_\bullet(F))$ the spectral radius of $H\!H_\bullet(F)$.
There is an isomorphism $H\!H_\bullet(X)\cong H^\bullet(X;\CC)$ given by 
the Hochschild--Kostant--Rosenberg isomorphism and the degeneration of the Hodge to de Rham spectral sequence.
Under this isomorphism, the endomorphism $H\!H_\bullet(\LL f^*)$ corresponds to 
the induced endomorphism $f^*$ of $H^\bullet(X;\CC)$ for each $f\in {\rm End}(X)$ 
(cf. Theorem~4.4, Lemma ~6.4 in \cite{Lunts} and references therein). 
In particular, we have $\rho(H\!H_\bullet(\LL f^*))=r(f)$.
Note also that $\rho([F])\le \rho(H\!H_\bullet(F))$ since $A_{num}^\bullet(X)_\CC:=A_{num}^\bullet(X)\otimes_\ZZ\CC$ 
is a quotient of a subspace $A_{hom}^\bullet(X)_\CC:=A_{hom}^\bullet(X)\otimes_\ZZ\CC$ 
of $\oplus_q H^{q,q}(X)\subset \oplus_{p,q} H^{p,q}(X)=H\!H_\bullet(X)$, where $A_{hom}^q(X)$ is 
the Chow group of codimension $q$ cycles on $X$ modulo (co)homological equivalence, 
and these linear maps are compatible with the induced actions of $F$. 

\subsection{Entropy and spectral radii}
From now on, we shall only consider the entropy $h_t(F)$ of an endofunctor $F$ of $\D^b(X)$ at $t=0$. 
For simplicity, set $\delta(G,F^nG'):=\delta_{0}(G,F^nG')$, $\delta'(G,F^nG'):=\delta'_{0}(G,F^nG')$, $h(F):=h_0(F)$ and so on. 
\begin{defn}
For a surjective endomorphism $f\in {\rm End}(X)$, the {\it categorical entropy} $h_{cat}(f)$ of $f$ is 
the entropy $h(\LL f^*)$ of its derived functor $\LL f^*:\D^b(X)\longrightarrow \D^b(X)$ at $t=0$. 
\end{defn}
\begin{rem}
Each surjective endomorphism $f\in {\rm End}(X)$, the induced endomorphism $H\!H_\bullet(\LL f^*)$ of $H\!H_\bullet(X)$ is invertible.
Therefore, we have $(\LL f^*)^nG\not\cong 0$ for all $n\ge 0$ and hence the categorical entropy $h_{cat}(f)$ is well-defined 
and $h_{cat}(f)\ge 0$. 
\end{rem}
In fact, it is possible to work in the category of smooth proper varieties over $\CC$, namely, 
for a smooth proper variety $Y$ over $\CC$ one can define the entropy $h(F)$ and the spectral radii 
$\rho([F])$, $\rho(H\!H_\bullet(F))$ of $F\in {\rm End}^{FM}(\D^b(Y))$ simply by setting $X=Y$
(however, we did not do this since our results are only for projective varieties). 
Concerning the entropy $h(F)$, there is a lower bound, under a certain technical assumption, by the natural logarithm of 
the spectral radius $\rho(H\!H_\bullet(F))$ of the induced $\CC$-linear endomorphism $H\!H_\bullet(F)$ of 
the Hochschild homology group $H\!H_\bullet(Y)$ (Theorem~2.9 in \cite{DHKK}). 
When $Y$ is projective, they also show that $h_{cat}(f)\le \rho(H\!H_\bullet(\LL f^*))$ 
(in fact, $h_{cat}(f)\le \rho([\LL f^*])$ holds, as we shall see later), 
which implies $h_{cat}(f)=\log r(f)$ for each surjective endomorphism $f\in {\rm End}(Y)$ satisfying the assumption (Theorem 2.12 in \cite{DHKK}).

From the above, it is natural to expect the following conjecture in general 
as a categorical generalization of the fundamental theorem on topological entropy due to Gromov-Yomdin (see Proposition~\ref{gy}). 
\begin{conj}\label{conj}
\footnote{
The lower bound $h(F)\geq\log\rho([F])$ in the conjecture does hold (\cite[eq.(1.1)]{KST}), 
but unfortunately, the upper bound $h(F)\leq\log\rho([F])$ does not in general (See \cite{Fan,Ouc}). 
Certainly, there is no effect on results in this paper. 
It is an interesting and important problem to find a characterization of functors attaining the lower bound.
}
Let $Y$ be a smooth proper variety over $\CC$.
For each $F\in {\rm End}^{FM}(\D^b(Y))$ such that the induced $\CC$-linear endomorphism 
$H\!H_\bullet(F)$ of $H\!H_\bullet(Y)$ is invertible, 
\begin{equation}
h(F)=\log\rho(H\!H_\bullet(F)).
\end{equation}
Moreover, if $Y$ is projective and $F$ is an autoequivalence, then 
\begin{equation}\label{nGY}
h(F)=\log\rho([F]).
\end{equation}
\end{conj}
\begin{rem}
The conjecture of this form does not hold in general for perfect derived categories of smooth proper dg algebras,  
due to phantom categories\footnote{The authors would like to thank A.I.~Efimov for helpful comments.}.
\end{rem}
\subsection{Results}
The following is the main theorem in this paper, which is a special case of Conjecture~\ref{conj}. 
\begin{thm}\label{main1}
Let $X$ be a smooth projective variety over $\CC$.
For each surjective endomorphism $f\in {\rm End}(X)$, 
\begin{eqnarray}
h_{cat}(f)&=&\lim_{n\to\infty}\frac{1}{n}\log|\chi(G,(\LL f^*)^nG^*)|\\
&=&\log\rho([\LL f^*])=\log\max_q r_q(f)=\log\max_q d_q(f).
\end{eqnarray}
In particular, we have $h_{cat}(f)=h_{top}(f)$.
\end{thm}
By Theorem~2.2 in \cite{Orl2}, any autoequivalence of $\D^b(X)$ is of Fourier--Mukai type.
If $d=1$, then the following holds.
\begin{prop}[Theorem~3.1 in \cite{Kik}]
Suppose that $X$ is a smooth projective curve.
Then, $h(F)=\log\rho([F])$ for all autoequivalences $F$ of $\D^b(X)$.
\end{prop}
This can be generalized to smooth projective varieties with the ample canonical or anti-canonical sheaf. 
\begin{thm}\label{main2}
Suppose that the canonical sheaf $\K_X$ or the anti-canonical sheaf $\K_X^{-1}$ is ample.
Then $h(F)=\log\rho([F])=0$ for all autoequivalences $F$ of $\D^b(X)$. 
\end{thm}
We shall give proofs of Theorems~\ref{main1} and \ref{main2} in the next section. 
\section{Proofs}
\subsection{Proof of Theorem~\ref{main1}}
First, note that the endofunctor $\LL f^*$ is of Fourier--Mukai type. 
For a very ample invertible sheaf $\L$, we set $G:=\oplus_{i=1}^{d+1}\L^i$ and $G^*:=\oplus_{i=1}^{d+1}\L^{-i}$. 
By Kodaira vanishing theorem, we have 
\[
\delta'(G,(\LL f^*)^nG^*)=(-1)^d\chi(G,(\LL f^*)^nG^*)=|\chi(G,(\LL f^*)^nG^*)|,
\]
which yields, by Proposition~\ref{smooth-proper}, 
\[
h_{cat}(f)=\lim_{n\to\infty}\frac{1}{n}\log\delta'(G,(\LL f^*)^nG^*)
=\lim_{n\to\infty}\frac{1}{n}\log|\chi(G,(\LL f^*)^nG^*)|.
\]
The Hirzebruch--Riemann--Roch theorem gives
\begin{equation*}
\chi(G,(\LL f^*)^nG^*)
=\int_X ch(G^*)ch((f^*)^n G^*) td(X).
\end{equation*}
It follows from explicit calculations that 
\begin{equation*}
(-1)^d\int_X ch(G^*)ch((f^*)^n G^*) td(X)
=\sum_{r=0}^{d}\sum_{q=0}^{d-r}c_{r,q}\int_X (f^*)^nc_1(\L)^q\cup c_1(\L)^{d-r-q}td_r(X),
\end{equation*}
where $c_{r,q}$ are rational numbers such that $c_{0,q}$ are all positive. 
Let $p$ be an integer such that $(r_p(f),l_p(f))=\displaystyle\max_{q} (r_q(f),l_q(f))$. 
We have
\begin{eqnarray*}
h_{cat}(f)&=&\lim_{n\to\infty}\frac{1}{n}\log|\chi(G,(\LL f^*)^nG^*)|\\
&=&\lim_{n\to\infty}\frac{1}{n}\log\left(
\sum_{r=0}^{d}\sum_{q=0}^{d-r}c_{r,q}\int_X (f^*)^nc_1(\L)^q\cup c_1(\L)^{d-r-q}td_r(X)\right)\\
&=&\log d_p(f)+\lim_{n\to\infty}\frac{1}{n}\log \left(
\sum_{r=0}^{d}\sum_{q=0}^{d-r} c_{r,q}\frac{\int_X (f^*)^nc_1(\L)^q \cup c_1(\L)^{d-r-q} td_r(X)}
{\int_X (f^*)^nc_1(\L)^p\cup c_1(\L)^{d-p}}\right). 
\end{eqnarray*}
Recall that $\int_X (f^*)^nc_1(\L)^q\cup c_1(\L)^{d-r-q} td_r(X)$ has at most the growth $n^{l_q(f)-1}r_q(f)^n$ 
as $n\rightarrow \infty$ (see Section~4, especially, \eqref{growth1}, \eqref{growth2} and Proposition~\ref{growth-p}).
Hence,
\begin{equation*}
\limsup_{n\to\infty}\left |\frac{\int_X (f^*)^nc_1(\L)^q \cup c_1(\L)^{d-r-q} td_r(X)}
{\int_X (f^*)^nc_1(\L)^p\cup c_1(\L)^{d-p}}\right |<\infty.
\end{equation*}
Note that $h_{cat}(f)$ and $\log d_p(f)$ do not depend on the choice of $\L$.
It follows that
\[
\lim_{n\to\infty}\frac{1}{n}\log \left(
\sum_{r=0}^{d}\sum_{q=0}^{d-r} c_{r,q}\frac{\int_X (f^*)^nc_1(\L)^q \cup c_1(\L)^{d-r-q} td_r(X)}
{\int_X (f^*)^nc_1(\L)^p\cup c_1(\L)^{d-p}}\right)=0, 
\]
since $c_{0,q}$ are all positive rational numbers and 
replacing $\L$ by sufficiently large power of $\L$ we may reduce the contributions from terms with $r>0$. 

Recall that $\N(X)_\CC$ can be identified with $A_{num}^\bullet(X)_\CC$. 
It follows that 
\[
\log\rho([\LL f^*])\le \log\displaystyle\max_q r_q(f)=\log r_p(f)
\]
since $A_{num}^\bullet(X)_\CC$ is a quotient of an $f^*$-invariant subspace $A_{hom}^\bullet(X)_\CC$ of $\oplus_{q}H^{q,q}(X)$ by an $f^*$-invariant subspace. 
Therefore, we have 
\[
\log\rho([\LL f^*])\le \log r_p(f)=\log d_p(f)=\lim_{n\to\infty}\frac{1}{n}\log|\chi(G,(\LL f^*)^nG^*)|=h_{cat}(f).
\]

On the other hand, it is easy to see that 
\[
h_{cat}(f)=\lim_{n\to\infty}\frac{1}{n}\log|\chi(G,(\LL f^*)^nG^*)|\leq\log\rho([\LL f^*])\leq\log r_p(f)=\log d_p(f).
\]
To summarize, we obtain $h_{cat}(f)=\log\rho([\LL f^*])=\log\displaystyle\max_q r_q(f)=\log\displaystyle\max_q d_q(f)$. 
\qed
\subsection{Proof of Theorem~\ref{main2}}
Let ${\rm Auteq}(\D^b(X))$ be the group of isomorphism classes of autoequivalences of $\D^b(X)$. 
It is shown by Bondal--Orlov that 
\[
{\rm Auteq}(\D^b(X))\cong ({\rm Aut}(X)\ltimes {\rm Pic}(X))\times\ZZ[1]
\]
 (Theorem~3.1 in \cite{BO}).
Choose an integer $m$ so that $\L:=\K_X^m$ is very ample. 
Set $G:=\oplus_{i=1}^{d+1}\L^i$ and $G^*:=\oplus_{i=1}^{d+1}\L^{-i}$.

For each automorphism $f\in{\rm Aut}(X)$, we have $\LL f^*G=G$ and hence 
\[
h_{cat}(f)=\lim_{n\to\infty}\frac{1}{n}\log\delta(G,(\LL f^*)^nG)=\lim_{n\to\infty}\frac{1}{n}\log\delta(G,G)=0.
\]

Each standard autoequivalence $F$ is represented as $F(-)=\LL f^*(-\otimes\L')[a]$ for an automorphism $f\in {\rm Aut}(X)$, an invertible sheaf $\L'$ and an integer $a$. 
Let $l$ be an integer such that $\L'':=\L'\otimes \K_X^l$ is anti-ample.  
Note that $(f^*)^n\L''$ is anti-ample for all $n\ge 0$ and the functors $(-\otimes\K_X^{-l})[a]$ and 
$\LL f^*(-\otimes\L'')$ commute. 
Since $h(-\otimes\K_X^{-l})=0$ by Lemma~\ref{pic}, we have $h(F)=h(\LL f^*(-\otimes\L''))$ by Lemma~\ref{entropy} (iv) and (v).  

Set $F'(-)=\LL f^*(-\otimes\L'')$. 
By Proposition~\ref{smooth-proper} and Kodaira vanishing theorem, 
\begin{eqnarray*}
h(F)&= & \lim_{n\to\infty}\frac{1}{n}\log\delta'(G,(F')^nG^*)\\
&=& \lim_{n\to\infty}\frac{1}{n}\log\delta'(G,G^*\otimes \L''\otimes f^*\L''\otimes\cdots\otimes(f^*)^{n-1}\L'')\\
&=& \lim_{n\to\infty}\frac{1}{n}\log|\chi(G,G^*\otimes \L''\otimes f^*\L''\otimes\cdots\otimes(f^*)^{n-1}\L'')|\\
&=& \lim_{n\to\infty}\frac{1}{n}\log|\chi(G,(F')^nG^*)|.
\end{eqnarray*}
It follows that 
\[
h(F)=\lim_{n\to\infty}\frac{1}{n}\log|\chi(G,(F')^nG^*)|\leq\log\rho([F'])=\log\rho([\LL f^*]).
\]
By Theorem~\ref{main1}, we have $\log\rho([\LL f^*])=h_{cat}(f)=0$ which yields the statement. 
\qed

\end{document}